\documentclass[twocolumn,amsthm]{autart}    
\usepackage{graphicx}
\usepackage{bm}
\usepackage{epsfig}
\usepackage{amssymb,amsfonts,amsmath}
\usepackage{subfigure}
\usepackage{setspace}
\usepackage{enumerate}
\usepackage{natbib}

\newcommand{\rd}{\mathrm{d}}
\bibpunct{(}{)}{;}{a}{,}{,}

\begin{document}

\begin{frontmatter}

\title{Ensemble Control of Stochastic Linear Systems\thanksref{footnoteinfo}} 

\author{Ji Qi}\ead{qij@ese.wustl.edu},    
\author{Anatoly Zlotnik}\ead{azlotnik@ese.wustl.edu},               
\author{Jr-Shin Li}\ead{jsli@seas.wustl.edu}  

\thanks[footnoteinfo]{Corresponding author Jr-Shin Li. Tel. 314-935-7340.
Fax 314-935-7500.}


\address{Electrical and Systems Engineering, Washington University, St. Louis, MO 63130}  


\begin{abstract}                          
	In this paper, we consider the problem of steering a family of independent, structurally identical, finite-dimensional stochastic linear systems with variation in system parameters between initial and target states of interest by using an open-loop control function. Our exploration of this class of control problems, which falls under the rising subject of ensemble control, is motivated by pulse design problems in quantum control. Here we extend the concept of ensemble control to stochastic systems with additive diffusion and jump processes, which we model using Brownian motion and Poisson counters, respectively, and consider optimal steering problems. We derive a Fredholm integral equation that is used to solve for the optimal control, which minimizes both the mean square error (MSE) and the error in the mean of the terminal state.  In addition, we present several example control problems for which optimal solutions are computed by numerically approximating the singular system of the associated Fredholm operator.  We use Monte Carlo simulations to illustrate the performance of the resulting controls. Our work has immediate practical applications to the control of dynamical systems with additive noise and parameter dispersion, and also makes an important contribution to stochastic control theory.

\end{abstract}

\end{frontmatter}

\section{Introduction} \label{sec_intro}
\label{sec:intro}
The behavior of physical, chemical, and biological systems can exhibit significant sensitivity to uncertainty or variation in system parameters.  This factor arises in practical control problems in many areas of science and engineering when there is uncertainty in the parameters of a single control system, or when a collection of structurally similar systems with variation in common parameters must be steered using a common control signal.  The former case has been approached by employing $H_\infty$ feedback control such as for quadratic stabilization of linear systems \citep*{xie92}, where parameter uncertainty is modeled with unknown time-dependent variations in the system dynamics.  When feedback is unavailable, however, analysis of these cases has given rise to the subject of ensemble control, which is motivated by the practical application of control theory to the fields of Nuclear Magnetic Resonance (NMR) spectroscopy and imaging (MRI) \citep*{li06thesis,li06pra}. 

Rapidly developing technologies based on quantum theory require the manipulation of large quantum ensembles, e.g., a spin ensemble on the order of Avogodro's number ($6 \times 10^{23}$), whose individual states cannot be measured and whose dynamics are subject to dispersion in natural frequencies. Such applications require open-loop controls whose effect is insensitive to dispersion in system parameters as well as to inhomogeneities in the applied radio-frequency (RF) control field \citep*{li06pra,li09tac}. A typical class of control problems significant to applications of NMR and MRI requires the design of RF pulses that drive a collection of spin systems with identical dynamics but parameter values unknown up to a given range between initial and target states, such that control performance is immune to physical parameter variation \citep*{kobzar05,pauly91}. A similar scenario appears in the field of neuroscience, where population level frequency control of neural oscillators requires open-loop stimuli that minimize signal power while achieving the desired response across a collection of neurons with variation in natural dynamics \citep*{brown04,li10,zlotnik11}.  Although initially motivated by the need to control large collections of similar systems, the paradigm of ensemble control can be readily used to solve any open-loop control problem in which the system response must be insensitive to uncertainty in model parameters.

Many practical control systems are subject to random dynamic disturbances, which are best treated from a probabilistic perspective by including a stochastic term in the system model.  We focus on the influence on ensemble control design of additive noise processes that perturb the system dynamics.  Because parameter variation can greatly affect the behavior of stochastic systems, especially when feedback is not used to attenuate disturbances, an analysis of ensemble control problems that include stochastic terms is crucial and meaningful.

It is important to clarify the distinction between parameter uncertainty and stochastic effects with respect to system dynamics.  Variation, or dispersion, in system parameters among members of a class of dynamical systems gives rise to differences in the actuation of individual members, each of which has fixed parameter values.  In practical applications, factors such as 
imprecise models or measurements and natural inhomogeneities result in uncertainty in parameters among a collection of systems, which is typically quantifiable within a given range.  In contrast, the state trajectory of a physical system often varies from a nominal path due to dynamic effects caused by the environment or inherent to the system itself, which may change randomly and depend on time.  These perturbations can be modeled as stochastic disturbances, whose nature may also depend on a parameter, so that the system state is described by a probability density function.  For a given parameter set, such stochastic effects can be attenuated by minimizing statistical objectives such as the expectation and variance of the system state distribution, as in the linear quadratic Gaussian (LQG) problem.   Our contribution is an alternative method for designing open-loop controls that are insensitive to parameter dispersion while optimizing these statistical objectives.

This paper is organized as follows. In the next section, we will provide a brief review of ensemble control for finite-dimensional time-varying deterministic linear systems, and introduce the definition of stochastic ensemble controllability. In Section \ref{sec_stoch}, we present our main result, which is a derivation of the conditions that the optimal ensemble control must satisfy in order to minimize the MSE at the terminal state. It is shown that the same control also minimizes the error in the mean, and coincides with the optimal control for the corresponding deterministic case as well. We provide a derivation for the case of Brownian diffusion, and indicate the necessary modifications for considering a Poisson type jump process. In Section \ref{sec_operator}, we discuss the numerical approximation of solutions to the Fredholm integral equation of the first kind that characterizes the optimal control by discretizing the kernel in time and parameter space, so that the singular system can be approximated by the singular value decomposition (SVD) of a matrix operator. Finally in Section \ref{sec_example}, several examples are given in order to illustrate the application and performance of our method. 

\section{Preliminaries} \label{sec_prelim}
\label{sec:imrt}
In this section, we review the basic ideas of ensemble control and introduce the notion of stochastic ensemble control.  Consider a parameterized family of finite-dimensional time-varying linear control systems
\begin{align}\label{ode}
\dot{X}(t,\beta)&=A(t,\beta)X(t,\beta)+B(t,\beta)u(t), \nonumber\\
   X & \in M \subset\mathbb{R}^n,  \,\, \beta\in K\subset\mathbb{R}^d, \,\, u\in U\subset \mathbb{R}^m,
\end{align}
where $A(t,\beta)\in\mathbb{R}^{n\times n}$ and $B(t,\beta)\in\mathbb{R}^{n\times m}$ have elements that are real $L_{\infty}$ and $L_2$ functions, respectively, defined on a compact set $D=[0,T]\times K\subset\mathbb{R}^2$, and are denoted $A\in L_\infty^{n\times n}(D)$ and $B\in L_2^{n\times m}(D)$.  The ensemble controllability conditions for the system (\ref{ode}) depend on the existence of an open-loop control $u:[0,T]\to U$ that can steer the instantaneous state of the ensemble $X:D\to M$ between any points of interest.  We consider the system (\ref{ode}) in a Hilbert space setting.  Let $\mathcal{H}_T=L_2^m[0,T]$ be the set of $m$-tuples, whose elements are vector-valued square-integrable functions defined on $0\leq t\leq T$, with an inner product defined by
\begin{align} \label{inpt}
\langle g,h \rangle_T = \int_0^T g'(t) h(t) \rd t,
\end{align}
where $'$ denotes the transpose.  Similarly, let $\mathcal{H}_K=L_2^n(K)$ be equipped with an inner product
\begin{align} \label{inpb}
\langle p,q \rangle_K = \int_K p'(\beta) q(\beta) \rd\mu(\beta),
\end{align}
where $\mu$ is the Lebesgue measure, and $K$ is the compact domain of the parameter $\beta$.  With well-defined addition and scalar multiplication, $\mathcal{H}_T$ and $\mathcal{H}_K$ are separable Hilbert spaces, where $||\cdot||_T$ and $||\cdot||_K$ denote their respective induced norms.  We now restate the definition of deterministic ensemble controllability \citep*{li11tac}.

\begin{defn}
	We say that the family (\ref{ode}) is ensemble controllable on the function space $\mathcal{H}_K$ if for all $\epsilon>0$, and all $X_0,\, X_F\in \mathcal{H}_K$, there exist $T>0$ and an open-loop piecewise-continuous control $u \in \mathcal{H}_T$, such that starting from $X(0,\beta)=X_0(\beta)$, the final state $X(T,\beta)=X_T(\beta)$ satisfies $||X_T-X_F||_K<\epsilon$.
\end{defn}

In other words, the system is ensemble controllable if it is possible to drive the system from $X_0$ to $X_F$ with respect to all $\beta\in K$, where the acceptable range of $T\in(0,\infty)$ may depend on $\epsilon$, $K$, and $U$.  Necessary and sufficient conditions have been provided for the ensemble controllability of finite-dimensional time-varying linear systems by applying integral operator theory \citep*{li11tac}. Given the initial state $X(0,\beta)=X_0(\beta)$ of the system (\ref{ode}), the variation of constants formula gives rise to 
\begin{align} \label{soldet}
X(T,\beta) & =   \Phi(T,0,\beta)X_0(\beta)  \nonumber \\ & \quad +\int_{0}^{T}\Phi(T,\sigma,\beta)B(\sigma,\beta)u(\sigma)\rd\sigma,
\end{align}
where $\Phi(T,0,\beta)$ is the transition matrix for the system $\dot{X}(t,\beta)=A(t,\beta)X(t,\beta)$.  The goal is for the terminal and target states to coincide in the function space $\mathcal{H}_K$, so setting $X(T,\beta)=X_F(\beta)$, pre-multiplying by $\Phi(0,T,\beta)$ and rearranging yields the integral operator equation
\begin{align} \label{critdet}
(Lu)(\beta) & \doteq\int_0^{T}\Phi(0,\sigma,\beta)B(\sigma,\beta)u(\sigma)\rd\sigma \nonumber \\ & =\Phi(0,T,\beta)X_F(\beta)-X_0(\beta),
\end{align}
which implicitly defines the solution $u$.  It is possible to use a spectral decomposition, called a singular system, of the operator $L$ to find an infinite eigenfunction series expansion for the function $u\in \mathcal{H}_T$ of minimum norm that satisfies (\ref{critdet}).  In practice, one can truncate the series and approximate the eigenfunctions numerically to estimate a solution that results in $||X_T-X_F||_K<\epsilon$.  The necessary and sufficient conditions for ensemble controllability are summarized and the method used to numerically approximate solutions to Fredholm operators is presented in detail in Section \ref{sec_operator}.

The controllability conditions for general, possibly nonlinear, ensemble control problems are presently unknown, and analytical control design methods remain a challenging problem, although analytical solutions exist for certain specific systems. An alternative direction is to use a powerful method for numerically solving optimal control problems called the pseudospectral method. This approach has been widely utilized to solve standard optimal control problems throughout the past decade for use in numerous applications \citep*{ross_legendre_2004}, and has recently been extended to solve problems in optimal ensemble control \citep*{Li_JCP11,Li_PNAS11}. The work presented here is focused on the derivation of analytical solutions, and numerical approximations, to optimal controls for  linear stochastic ensemble systems.


To investigate the effect of stochastic factors on the controllability of ensembles, we extend the notion of ensemble controllability to incorporate systems where the state is a random variable on a given probability space.   Consider a modification to the collection (\ref{ode}) by defining a family of stochastic control systems
\begin{align}\label{odes}
\rd X(t,\beta,\omega) &=A(t,\beta)X(t,\beta,\omega) \rd t \nonumber \\ & \quad +B(t,\beta)u(t)dt+G(t,\beta) \rd S(t,\omega),
\end{align}
for $X\in M \subset\mathbb{R}^n$, $\beta\in K\subset\mathbb{R}^d$, $u\in U\subset \mathbb{R}^m$, and $S\in V \subset\mathbb{R}^k$,  where $A(t,\beta)\in\mathbb{R}^{n\times n}$, $B(t,\beta)\in\mathbb{R}^{n\times m}$, and $G(t,\beta)\in\mathbb{R}^{n\times k}$ have elements that are real $L_{\infty}$, $L_2$, and $L_2$ functions, respectively, defined on a compact set $D=[0,T]\times K\subset\mathbb{R}^2$, and are denoted $A\in L_\infty^{n\times n}(D)$, $B\in L_2^{n\times m}(D)$, and $G\in L_2^{n\times k}(D)$.  The equation (\ref{odes}) is stated in the It\^o sense, where the term $\rd S$ is the differential of a continuous-time stochastic process $S:[0,T]\times \Omega\to V \subset \mathbb{R}^k$ on a probability space $(\Omega,\mathcal{F},\mathbb{P})$. The ensemble controllability conditions for the systems (\ref{odes}) depend on the existence of an open-loop control $u\in \mathcal{H}_T$ that can steer the instantaneous expected state of the stochastic ensemble $\mathcal{E} X:D\to M$ between any points of interest in the function space $\mathcal{H}_K$, where $\mathcal{E}$ denotes expectation over the stochastic space $\Omega$.  Our analysis of stochastic ensemble controllability of the collection (\ref{odes}) is based on the approach taken in the deterministic case, hence a necessary condition is the ensemble controllability of the corresponding deterministic system, i.e. when $G$ is identically zero.  This requirement gives rise to the following definition.

\begin{defn} \label{stochcontr}
We say that the family of systems (\ref{odes}) is ensemble controllable on the function space $\mathcal{H}_K$ if for all $\epsilon>0$, and all $X_0, X_F\in \mathcal{H}_K$, there exist $T>0$ and an open-loop piecewise-continuous control $u\in\mathcal{H}_T$ such that starting from $X(0,\beta)=X_0(\beta)$, the final state $X(T,\beta)=X_T(\beta)$ satisfies $||\mathcal{E}X_T-X_F||_K<\epsilon$.
\end{defn}

If the system (\ref{odes}) is ensemble controllable, then a natural objective for the problem of steering the ensemble from $X_0$ to $X_F$ in time $T$ is 
to minimize the error in the mean of the terminal state with respect to the target state,
\begin{align} \label{merror}
\mathcal{J}_1 &= ||\mathcal{E}X_T-X_F||_K \nonumber \\ & = \left(\int_K||X_T(\beta)-X_F(\beta)||_2^2\,\rd\mu(\beta)\right)^{1/2},
\end{align}
where $||\cdot||_2$ denotes the Euclidean norm.  In addition to error in the mean, another important factor in stochastic control is MSE, which is minimized using the objective
\begin{align} \label{mse}
\mathcal{J}_2 & =  \mathcal{E}||X_T-X_F||_K^2 \\ & = \int_K\mathcal{E}||X(T,\beta)-X_F||_2^2 \rd\mu(\beta), \nonumber
\end{align}
where we invoke Fubini's theorem to evaluate the expectation before integrating with respect to $\beta$.
\begin{rem}  If the $A$, $B$, $G$, $X_0$, and $X_F$ all depend continuously on the parameter $\beta$, then the $L_\infty$ norm can be used in Definition \ref{stochcontr} and the objectives $\mathcal{J}_1$ and $\mathcal{J}_2$. In that case, the objective $\mathcal{J}_1$ is endowed the physical meaning of guaranteeing that the error in the mean is bounded irrespective of $\beta\in K$, and the objective $\mathcal{J}_2$ guarantees that the MSE is uniformly bounded with respect to $\beta$.  In order to broaden the applicability of our approach, we have chosen to use the $L_2$ norm, at the cost of uniformity in this respect.
\end{rem}
In the following section, we will prove that the optimal controls that minimize objectives (\ref{merror}) and (\ref{mse}) coincide, and also solve the deterministic control problem when $S$ is a Brownian diffusion process.

\section{Stochastic Ensemble Control} \label{sec_stoch}

In this section, we derive optimal ensemble controls for systems with additive diffusion and jump processes, which we model using Brownian motion and Poisson counters, respectively.  The solution in each case is obtained implicitly in the form of a Fredholm integral equation of the first kind, so it is possible for specific cases to obtain an analytical expression of the control function as an infinite series of weighted eigenfunctions of the Fredholm operator as in the deterministic case \citep*{li11tac}.  We present the derivation for the case of Brownian diffusion first, and discuss the modifications necessary for considering a Poisson type jump process.

\subsection{Brownian Diffusion}

Consider a finite-dimensional time-varying stochastic linear system with an additive Brownian noise process governed by the It\^o differential equation
\begin{align}\label{sysbm}
\rd X(t,\beta,\omega) &=A(t,\beta)X(t,\beta,\omega)\rd t \nonumber \\ & \quad +B(t,\beta)u(t) \rd t+G(t,\beta) \rd W(t,\omega),
\end{align}
for $ X\in M  \subset  \mathbb{R}^n$, $\beta\in K\subset\mathbb{R}^d$, $u\in U\subset \mathbb{R}^m$, and $W\in \mathbb{R}^k$, where $A(t,\beta)\in\mathbb{R}^{n\times n}$, $B(t,\beta)\in\mathbb{R}^{n\times m}$, and $G(t,\beta)\in\mathbb{R}^{n\times k}$ have elements that are real $L_{\infty}$, $L_2$, and $L_2$ functions, respectively, defined on a compact set $D=[0,T]\times K\subset\mathbb{R}^2$, and are denoted $A\in L_\infty^{n\times n}(D)$, $B\in L_2^{n\times m}(D)$, and $G\in L_2^{n\times k}(D)$.  The vector valued stochastic process $W:[0,T]\times \Omega \to \mathbb{R}^k$ consists of independent identically distributed (i.i.d.) standard normal random variables on the probability space $(\Omega,\mathcal{F},\mathbb{P})$, and has the natural filtration so that $\mathcal{E}W(t)=0$ and $\mathcal{E}[W(t)W'(s)]=I\min(t,s)$.  We use $\mathcal{E}$ to denote expectation over the space $\Omega$ with respect to measure $\mathbb{P}$, and we omit $\omega$ as an argument of $X$ and $W$ from now on in order to adhere to a straightforward notation for stochastic calculus \citep*{brockett02}.  We first derive the control that minimizes $\mathcal{J}_1$ in (\ref{merror}) as follows.

\begin{prop} \label{prop1bm} 
	For ensemble controllable system (\ref{sysbm}), the control $u\in \mathcal{H}_T$ that minimizes $\mathcal{J}_1$, the norm of the error in the mean of the terminal state (\ref{merror}), while steering the system from the initial state  $X(0,\beta)=X_0(\beta)$ to the target state $X(T,\beta)=X_F(\beta)$ must satisfy (\ref{critdet}).
\end{prop}

\begin{pf} Using Fubini's theorem and the fact $\mathcal{E}W(t)=0$,
\begin{align}
\mathcal{E}\left[\int_0^TG(t,\beta)dW(t)\right] = \int_0^TG(t,\beta)\mathcal{E}\rd W(t) = 0,
\end{align}
hence the expectation of the terminal state is given by
\begin{align} \label{solbm}
\mathcal{E}X(T,\beta) &=\Phi(T,0,\beta)X_0(\beta) \nonumber \\ & \quad +\int_{0}^{T}\Phi(T,\sigma,\beta)B(\sigma,\beta)u(\sigma)d\sigma,
\end{align}
which coincides with the solution (\ref{soldet}) in the deterministic case.  Setting $\mathcal{E}X(T,\beta)=X_F(\beta)$ results in (\ref{critdet}). 
\end{pf}
This proposition implies that if (\ref{sysbm}) is ensemble controllable then one can achieve $\mathcal{J}_1<\epsilon$, so that the mean of the terminal state approximates the target in the space $\mathcal{H}_K$.  The minimization of the MSE of the terminal state, given by objective (\ref{mse}), is addressed in the following theorem, which is the main theoretical result of this paper.

\begin{thm} \label{thmain} 
	For ensemble controllable system (\ref{sysbm}), the control $u\in \mathcal{H}_T$ that steers the system from the initial state $X(0,\beta)=X_0(\beta)$ to the terminal state $X(T,\beta)=X_F(\beta)$, while minimizing $\mathcal{J}_2$, the mean square error in the terminal state (\ref{mse}), must satisfy (\ref{critdet}) for all $\beta\in K$. The minimum value of $\mathcal{J}_2$ is $\int_K\mathrm{tr}\,C(T,\beta)\rd\mu(\beta)$,	where $C(t,\beta)=\int_0^t\Phi(t,\sigma,\beta) G(\sigma,\beta)G'(\sigma,\beta)\Phi'(t,\sigma,\beta)\rd\sigma$ and $\Phi(t,0,\beta)$ is the transition matrix for $\dot{X}(t,\beta)=A(t,\beta)X(t,\beta)$.
\end{thm}


\begin{pf}
Let $\Phi(t,0,\beta)$ denote the transition matrix of $A(t,\beta)$.  The mean square error in the terminal state $X_T(\beta)=X(T,\beta)$ as a function of $\beta$ can be expressed as
\begin{align} \label{mseext}
\mathcal{E}||X_T-X_F||_2^2=\mathcal{E}X_T'X_T-2X_F'\mathcal{E}X_T+X_F'X_F,
\end{align}
where we omit the argument $\beta$ in order to simplify notation, and $'$ denotes the transpose.  Noting that $\mathcal{E}X'X=\mathrm{tr}\mathcal{E}XX'$, where $\mathrm{tr}$ denotes the trace of a square matrix, we set
$Z=XX'$ and apply It\^o's rule to obtain
\begin{align}
\rd Z &=(AZ +Z A'+BuX'+Xu'B'+GG')\rd t \nonumber \\ &\quad +(GX'+XG')\rd W,
\end{align}
and taking the expectation results in the deterministic linear matrix equation
\begin{align} \label{sigmaeqbm}
\dot{\Sigma}=A\Sigma+\Sigma A'+Bu(\mathcal{E}X)'+\mathcal{E}Xu'B'+GG',
\end{align}
where $\Sigma = \mathcal{E}Z=\mathcal{E}XX'$. Solving (\ref{sigmaeqbm}) gives rise to
\begin{align} \label{sigmasolbm}
\Sigma(t)&=\Phi(t,0)\Sigma(0)\Phi'(t,0) \nonumber \\
&  +\int_0^t\Phi(t,\sigma) B(\sigma)u(\sigma)\mathcal{E}X'(\sigma)\Phi'(t,\sigma)\rd\sigma\nonumber \\
&  +\int_0^t\Phi(t,\sigma)\mathcal{E}X(\sigma)u'(\sigma)B'(\sigma)\Phi'(t,\sigma)\rd\sigma \nonumber \\
&  +\int_0^t\Phi(t,\sigma)G(\sigma)G'(\sigma)\Phi'(t,\sigma)\rd\sigma
\end{align}
as a result of elementary linear system theory \citep*{brockett70}.  Observe that the quantities
\begin{align} \label{indebm}
C(t)=& \int_0^t\Phi(t,\sigma) G(\sigma)G'(\sigma)\Phi'(t,\sigma)\rd\sigma
\end{align}
and $\Phi(t,0)\Sigma(0)\Phi'(t,0)$ are independent of the control $u$.  Let $\Phi_T$ denote $\Phi(T,0)$, and define a change of variables
\begin{align} \label{ydef}
Y(\sigma)=\int_0^{\sigma}\Phi(0,\tau)B(\tau)u(\tau)\rd\tau,
\end{align}
so that $Y(0)=0$. Substituting (\ref{solbm}) for $\mathcal{E}X(\sigma)$ in equation (\ref{sigmasolbm}) leads to
\begin{align} \label{sigmasol2bm}
\Sigma(T) &= \Phi_T\int_0^T\dot{Y}(\sigma)(X_0+Y(\sigma))'\rd\sigma \, \Phi_T'
\nonumber \\
& +\Phi_T\int_0^T(X_0+Y(\sigma))\dot{Y}'(\sigma)\rd\sigma \, \Phi_T' \nonumber \\
&   + \Phi_T\Sigma(0)\Phi'_T  +C(T).
\end{align}
This, together with $Y(0)=0$, yields
\begin{align} \label{msebm}
\mathcal{E}X_T'X_T &=  \mathrm{tr}\,\int_0^T(X_0+Y(\sigma))'\Phi_T'\Phi_T\dot{Y}(\sigma)\rd\sigma
\nonumber \\
& \quad + \mathrm{tr}\,\int_0^T\dot{Y}'(\sigma)\Phi_T'\Phi_T(X_0+Y(\sigma))\rd\sigma \nonumber \\
& \quad + \mathrm{tr}\,\Phi_T\Sigma(0)\Phi'_T  + \mathrm{tr}\,C(T) \nonumber\\
 &=2X_0'\Phi_T'\Phi_T(Y(T)-Y(0)) \nonumber \\
& \quad  +Y'(T)\Phi_T'\Phi_TY(T)  -Y'(0)\Phi_T'\Phi_TY(0) \nonumber \\
& \quad  +X_0'\Phi'_T\Phi_TX_0+\mathrm{tr}\,C(T) \nonumber \\
 &=2X_0'\Phi_T'\Phi_TY(T)+Y'(T)\Phi_T'\Phi_TY(T) \nonumber \\
& \quad  +X_0'\Phi'_T\Phi_TX_0+\mathrm{tr}\,C(T) .
\end{align}
Because $X_0$ is non-random, $\mathrm{tr} [\Phi_T \Sigma(0) \Phi_T] = ||\Phi_TX_0||_2^2$, so we can compute the mean square error by substituting (\ref{msebm}) and (\ref{solbm}) into (\ref{mseext}) to obtain
\begin{align} \label{msebm2}
\mathcal{E} ||&X_T-X_F||_2^2 \nonumber \\
&=  \, \mathrm{tr} \, (\Sigma(T))-2X_F'\Phi_T(X_0+Y(T)) +X_F'X_F \nonumber\\ &=  \, 2X_0'\Phi_T'\Phi_TY(T)+Y'(T)\Phi_T'\Phi_TY(T)  \nonumber \\
& \quad -2X_F'\Phi_T(X_0+Y(T))+X_F'X_F \nonumber \\
& \quad +X_0'\Phi'_T\Phi_TX_0+\mathrm{tr}\,(C(T))  \nonumber \\
&=  \, ||\Phi_TY(T)+\Phi_TX_0-X_F||_2^2  +\mathrm{tr}\,C(T).
\end{align}
Equation (\ref{msebm2}) attains its minimum with respect to $u$ when
$Y(T)+X_0-\Phi_T^{-1}X_F=0$, that is, when
\begin{align} \label{msecrit} 
\int_0^{T}\Phi(0,\tau,\beta)&B(\tau,\beta)u(\tau)\rd\tau  \nonumber \\
&= \Phi(0,T,\beta)X_F(\beta)-X_0(\beta),
\end{align}
which is the condition (\ref{critdet}).  When this is satisfied for all $\beta\in K$, then the objective $\mathcal{J}_2$ is minimized at the minimum value
\begin{align} \label{obj2min}
\min_u\mathcal{J}_2=\int_K\mathrm{tr}\,C(T,\beta)\rd\mu(\beta).
\end{align}
\end{pf}

\begin{rem} \label{remkbm} Theorem \ref{thmain} together with Proposition \ref{prop1bm} demonstrate that the control minimizing the error in the mean of the terminal state also minimizes the mean square error in the terminal state for all systems in the ensemble (\ref{sysbm}).  The triangle inequality yields
$\mathcal{E}||X_T-X_F||_K \leq \mathcal{E} ||X_T-\mathcal{E} X_T||_K + ||\mathcal{E} X_T-X_F||_K$,
where equality is achieved when $\mathcal{E} X_T-X_F = 0$ in the function space $\mathcal{H}_K$.  In that case $\mathcal{E}||X_T-X_F||_K= \mathcal{E} ||X_T-\mathcal{E} X_T||_K$.  It follows that minimizing $\mathcal{J}_1$ also minimizes $\mathcal{J}_2$.
\end{rem}

\subsection{Poisson Jump Process}

Consider a finite-dimensional time-varying stochastic linear system with an additive Poisson jump process governed by the It\^o differential equation
\begin{align}\label{syspp}
\rd X(t,\beta,\xi) &=A(t,\beta)X(t,\beta,\omega)\rd t+B(t,\beta)u(t) \rd t \nonumber \\ & \quad +G(t,\beta)\rd N(t,\xi)
\end{align}
for $X\in M \subset\mathbb{R}^n$, $\beta\in K\subset\mathbb{R}$, $u\in U\subset \mathbb{R}^m$, and $W\in \mathbb{R}^k$, where $A(t,\beta)\in\mathbb{R}^{n\times n}$, $B(t,\beta)\in\mathbb{R}^{n\times m}$, and $G(t,\beta)\in\mathbb{R}^{n\times k}$ have elements that are real $L_{\infty}$, $L_2$, and $L_2$ functions, respectively, defined on a compact set $D=[0,T]\times K\subset\mathbb{R}^2$, and are denoted $A\in L_\infty^{n\times n}(D)$, $B\in L_2^{n\times m}(D)$, and $G\in L_2^{n\times k}(D)$.  The vector valued stochastic process $N:[0,T]\times \Xi \to \mathbb{R}^k$ consists of i.i.d. Poisson counters on the appropriate probability space $(\Xi,\mathcal{G},\mathbb{N})$, with a vector $\lambda\in\mathbb{R}^k$ of constant intensities, and with the natural filtration so that $\mathcal{E}N(t)=\lambda t$ and $\mathrm{Cov}(N(t),N(s))=\Lambda\min(t,s)$.  Here $\Lambda$ denotes the diagonal matrix containing the elements of $\lambda$ in its main diagonal, and zero elsewhere.  We use $\mathcal{E}$ to denote expectation over the space $\Xi$ with respect to measure $\mathbb{N}$, and we omit $\xi$ as an argument of $X$ and $N$ from now on in order to adhere to a straightforward notation for stochastic calculus \citep*{brockett02}.  We first derive the control that minimizes $\mathcal{J}_1$ in (\ref{merror}) as follows.

\begin{prop} \label{prop1pp} 
	Given the ensemble system (\ref{syspp}), the control $u\in \mathcal{H}_T$ which minimizes $\mathcal{J}_1$, the error in the mean of the terminal state (\ref{merror}), while steering the system from the initial state $X(0,\beta)=X_0(\beta)$ to the target state $X(T,\beta)=X_F(\beta)$ must satisfy
\begin{align} \label{critpp}
(Ru)(\beta) &\doteq \int_0^{T}\Phi(0,\sigma,\beta)(B(\sigma,\beta)u(\sigma)+G(\sigma,\beta)\lambda)\rd\sigma \nonumber \\ &
=\Phi(0,T,\beta)X_F(\beta)-X_0(\beta).
\end{align}
\end{prop}
\begin{pf} The property $\mathcal{E}N(t)=\lambda t$ and Fubini's theorem result in
\begin{align}
\mathcal{E}\left[\int_0^TG(t,\beta)\rd N(t)\right] &= \int_0^TG(t,\beta)\lambda \rd t,
\end{align}
hence the expectation of the terminal state is given by
\begin{align} \label{solpp}
\mathcal{E}X(T,\beta)&=\Phi(T,0,\beta)X_0(\beta) \nonumber \\ & \quad +\int_{0}^{T}\Phi(T,\sigma,\beta)B(\sigma,\beta)u(\sigma)\rd\sigma \nonumber \\ &  \quad
+\int_{0}^{T}\Phi(T,\sigma,\beta)G(\sigma,\beta)\lambda\rd\sigma.
\end{align}
Setting $X(T,\beta)=X_F(\beta)$, pre-multiplying by $\Phi(0,T,\beta)$ and rearranging yields integral operator equation (\ref{critpp}). 
\end{pf}
This proposition implies that if the family (\ref{sysbm}) is ensemble controllable then it is possible to achieve $\mathcal{J}_1<\epsilon$, so that the mean of the terminal state approximates the target state in the space $\mathcal{H}_K$.  The minimization of the mean square error of the terminal state, given by objective (\ref{mse}), is addressed in the following theorem, which is analogous to the main result.

\begin{thm} \label{thpp} 
	Given the ensemble controllable system (\ref{syspp}), the control $u\in \mathcal{H}_T$ which steers the system from the initial state  $X(0,\beta)=X_0(\beta)$ to the terminal state $X(T,\beta)=X_F(\beta)$, while minimizing $\mathcal{J}_2$, the mean square error in the terminal state (\ref{mse}), must satisfy (\ref{critpp}) for all $\beta\in K$. The minimum value of $\mathcal{J}_2$ is $\int_K\mathrm{tr}\,C(T,\beta)\rd\mu(\beta)$,	where $C(t,\beta)=\int_0^t\Phi(t,\sigma,\beta) G(\sigma,\beta)\Lambda G'(\sigma,\beta)\Phi'(t,\sigma,\beta)\rd\sigma$ with $\Lambda=\mathrm{diag} (\lambda)$ being the diagonal matrix of intensities, and $\Phi(t,0,\beta)$ is the transition matrix for $\dot{X}(t,\beta)=A(t,\beta)X(t,\beta)$.
\end{thm}


\begin{pf}
We can rewrite (\ref{syspp}) in the form $\rd X =AX \rd t+Bu\rd t+\sum_{i=1}^kg_i\rd N_i$,  where $g_i$ is the $i^{\mathrm{th}}$ column of $G$, and $\rd N_i$ is the $i^{\mathrm{th}}$ element of $\rd N$. Setting $Z=XX'$ and applying It\^o's rule for jump processes gives rise to
\begin{align}
dZ &=(AZ+Z A'+BuX'+Xu'B')\rd t  \nonumber \\ & \quad + \sum_{i=1}^k(g_iX'+Xg_i'+g_ig_i')\rd N_i \nonumber\\
&=(AZ+Z A'+BuX'+Xu'B')\rd t \nonumber \\ & \quad + G\rd N X'+X\rd N'G'+G(\mathrm{diag} (\rd N)) G',
\end{align}
where $\mathrm{diag} (\rd N) \in\mathbb{R}^{k\times k}$ is a diagonal matrix with $\rd N$ on the main diagonal. Taking the expectation results in the deterministic linear matrix equation
\begin{align} \label{sigmaeqpp}
\dot{\Sigma} &=A\Sigma+\Sigma A'+[Bu+ G \lambda](\mathcal{E}X)'  \nonumber \\ & \quad +\mathcal{E}X[Bu+ G\lambda]'+G\Lambda G',
\end{align}
where $\Sigma = \mathcal{E}XX'$ and $\Lambda=\mathrm{diag} (\lambda)$ is the diagonal matrix of intensities. By applying the matrix variation of constants formula \citep*{brockett70}, we obtain
\begin{align} \label{sigmasolpp}
\Sigma(t)&=\Phi(t,0)\Sigma(0)\Phi'(t,0) \nonumber \\
& +\int_0^t\Phi(t,\sigma) (B(\sigma)u(\sigma) + G(\sigma)\lambda)\mathcal{E}X'(\sigma)\Phi'(t,\sigma)\rd\sigma \nonumber \\ &  +\int_0^t\Phi(t,\sigma)\mathcal{E}X(\sigma)(B(\sigma)u(\sigma) + G(\sigma)\lambda)'\Phi'(t,\sigma)\rd\sigma \nonumber \\ &
+ \int_0^t\Phi(t,\sigma)G(\sigma)\Lambda G'(\sigma)\Phi'(t,\sigma)\rd\sigma.
\end{align}
We denote the quantity
\begin{align} \label{indepp}
C(t)=& \int_0^t\Phi(t,\sigma) G(\sigma)\Lambda G'(\sigma)\Phi'(t,\sigma)\rd\sigma
\end{align}
and apply a change of variables
\begin{align} \label{ydefpp}
Y(\sigma)=\int_0^{\sigma}\Phi(0,\tau)[B(\tau)u(\tau)+G(\tau)\lambda]\rd\tau,
\end{align}
so that substituting (\ref{solpp}) for $\mathcal{E}X(T,\beta)$ in (\ref{sigmasolpp}) leads to
(\ref{sigmasol2bm}), as in the case of Brownian diffusion, and the proof continues using the steps in the proof of Theorem \ref{thmain}. It can be shown that $\mathcal{J}_2$ attains its minimum with respect to $u$ when
$Y(T)+X_0-\Phi_T^{-1}X_F=0$, that is, when
\begin{align} \label{msecritpp} 
\int_0^{T}\Phi(0,\tau,\beta)&[B(\tau,\beta)u(\tau)+G(\tau,\beta)\lambda]\rd\tau \nonumber \\ & = \Phi(0,T,\beta)X_F(\beta)-X_0(\beta),
\end{align}
which is (\ref{critpp}).  When this holds for all $\beta\in K$, the objective $\mathcal{J}_2$ is minimized at value (\ref{obj2min}) with $C$ given in (\ref{indepp}). 
\end{pf}

\begin{rem} The proof of Theorem \ref{thpp} together with Proposition \ref{prop1pp} demonstrate that the control minimizing the error in the mean of the terminal state also minimizes the mean squared error in the terminal state for all systems in the ensemble (\ref{syspp}).
\end{rem}

The above theorems form a theoretical foundation for ensemble control of stochastic linear systems.  In the following section, we describe a numerical method for solving the equations (\ref{critdet}) and (\ref{critpp}) that enables the synthesis of optimal ensemble controls. 

\section{Ensemble Control Synthesis} \label{sec_operator}

We characterize ensemble controllability of a family of finite-dimensional time-varying linear systems by certain conditions on the associated input-to-state operator \citep*{li11tac}.  The system of interest (\ref{ode}) is considered in the Hilbert space setting described in Section \ref{sec_prelim}.  The Hilbert spaces $\mathcal{H}_T$ and $\mathcal{H}_K$ with inner products defined by (\ref{inpt}) and (\ref{inpb}), respectively, have well-defined addition and scalar multiplication operations, and are thus both separable.  It follows that  the operator defined in equation (\ref{critdet}) satisfies $L\in\mathcal{B}(\mathcal{H}_T,\mathcal{H}_K)$, where $\mathcal{B}(\mathcal{H}_T,\mathcal{H}_K)$ is the set of bounded linear operators from $\mathcal{H}_T$ to $\mathcal{H}_K$.  As a consequence, $L$ has the adjoint $L^*$ satisfying
\begin{align}
\langle \xi, Lu \rangle_K=\langle L^*\xi, u \rangle_T, \quad \forall \,\, \xi\in \mathcal{H}_K, \quad u\in \mathcal{H}_T.
\end{align}
We omit the subscript on the inner product when it is clear from the context.  The conditions characterizing ensemble controllability are related to the singular system of the operator $L$, which requires the following important definition.

\begin{defn} \label{singsys} {\it Singular System} \citep*{gohberg03}: Let $Y$ and $Z$ be Hilbert spaces and $K:Y\to Z$ be a compact operator.  If $(\sigma_n^2,\nu_n)$ is an eigensystem of $KK^*$ and $(\sigma_n^2,\mu_n)$ is an eigensystem of $K^*K$, namely, $KK^*\nu_n=\sigma_n^2\nu_n$, $\nu_n\in Z$, and $K^*K\mu_n=\sigma_n^2\mu_n$, $\mu_n\in Y$, where $\sigma_n>0$, and the two systems are related by
\begin{align} \label{sscond}
K\mu_n=\sigma_n\nu_n, & \quad K^*\nu_n=\sigma_n\mu_n,
\end{align}
we say that $(\sigma_n,\mu_n,\nu_n)$ is a singular system of $K$.
\end{defn}

\begin{rem} It can be shown that the operator $L$ in (\ref{critdet}) is compact \citep*{li11tac}, and hence $LL^*$ and $L^*L$ are both compact, self-adjoint, and nonnegative operators.  By the Spectral theorem, $L^*L$ can be represented in terms of its positive eigenvalues, namely, $L^*Ly=\sum_n\sigma_n^2\langle y,\mu_n \rangle \mu_n$ for all $y\in \mathcal{H}_T$.  Moreover, because $L^*L\mu_n=\sigma_n^2\mu_n$, the relations $L\mu_n=\sigma_n\nu_n$ and $L^*\nu_n=\sigma_n\mu_n$ follow by taking $\nu_n=(1/\sigma_n)L\mu_n$.  This is the infinite dimensional analogue of the singular value decomposition of a matrix.
\end{rem}

The following theorem provides necessary and sufficient conditions for ensemble controllability of finite-dimensional time-varying linear systems with input-to-state operator $L$, and is valid for stochastic systems when the appropriate compact operator is interposed.

\begin{thm} \textup{\citep*{li11tac}} \label{thmaindet}
The family of systems (\ref{ode}) is ensemble controllable on the function space $\mathcal{H}_K$ if and only if for any given initial and final state $X_0 \, X_F\in \mathcal{H}_K$, at time $t=0$ and $t=T<\infty$ respectively, and for $\xi=\Phi(0,T,\beta)X_F-X_0$, the conditions
\begin{equation}\begin{array}{rl}
 (i) \displaystyle \quad\sum_{n=1}^{\infty}\frac{|\langle\xi,\nu_n\rangle|^2}{\sigma_n^2}<\infty, \quad
(ii) \displaystyle \quad \xi\in\overline{\mathcal{R}(L)},\end{array}
\end{equation}
hold, where $L$ is the input-to-state operator of the system (\ref{ode}) defined in equation (\ref{critdet}),
with $\overline{\mathcal{R}(L)}$ denoting the closure of the range space of $L$, and the collection of triples $(\sigma_n,\mu_n,\nu_n)$ is a singular system of the linear operator $L$.  Moreover, the control law
\begin{align} \label{seqmaindet}
u=\sum_{n=1}^\infty \frac{1}{\sigma_n}\langle \xi, \nu_n\rangle\mu_n
\end{align}
satisfies $\langle u,u \rangle \leq \langle u_0,u_0 \rangle$ for all $u_0\in \mathcal{U}$ and $u_0\neq u$, where $\mathcal{U}=\{v\,|\,Lv=\xi \text{ with } (i) \text{ and } (ii)\}$.
\end{thm}

The above theorem, together with Definition \ref{singsys}, enables a method for approximating solutions to (\ref{critdet}) of minimum norm \citep*{zlotnik12acc}.  The integral operator equation $Lu=\xi$ can be approximated by a linear matrix equation, so that the singular system $(\sigma_n,\mu_n,\nu_n)$ of $L$ can be approximated by the singular value decomposition of this matrix.  Let $\{\beta_j\}$ be a collection of points that uniformly distributed throughout the space $K$ for $j=0,1,2,\ldots,P$, and let $\{t_k\}$ be a collection of points that linearly interpolate the time domain $[0,T]$ for $k=0,1,\ldots,N$, including endpoints, with $t_{k}-t_{k-1}=\delta$.  Using this grid of nodes, we approximate
\begin{align}
(Lg)(\beta)&=\int_0^T\Phi(0,t,\beta)B(t,\beta)g(t)\rd t \nonumber \\
& = \sum_{k=1}^N \int_{t_{k-1}}^{t_k}\Phi(0,t,\beta)B(t,\beta)g(t) \rd t \nonumber\\
& \approx \sum_{k=1}^N\delta\Phi(0,t_k,\beta)B(t_k,\beta)g(t_k)
\end{align}
for each $\beta\in\{\beta_j\}$. Hence the action of the operator $L$ on a function $g\in \mathcal{H}_T$ can be approximated by the action of a block matrix $W\in\mathbb{R}^{nP\times mN}$, with $n\times m$ blocks $W_{jk}=\delta \Phi(0,t_k,\beta_j)B(t_k,\beta_j)$, on a vector $\hat{g}\in \mathbb{R}^{mN}$, with $N$ blocks   $\hat{g}_{k}=g(t_k)$ of dimension $m\times 1$.  If the SVD of this matrix is $W=U\Sigma V'$, and $\bar{u}_j$ and $\bar{v}_j$ are columns of $U$ and $V$, respectively, corresponding to the singular value $s_j$, then $WW' \bar{u}_j=s_j^2\bar{u}_j$ and $W' W\bar{v}_k=s_j^2\bar{v}_k$.  Therefore the SVD $(s_j,\bar{v}_j,\bar{u}_j)$ of the matrix $W$ approximates the singular system $(\sigma_j,\mu_j,\nu_j)$ of the operator $L$, where $\bar{v}_j$ and $\bar{u}_j$ are discretizations of $\mu_j$ and $\nu_j$, respectively.  Now suppose that $\hat{\xi}\in\mathbb{R}^{nP}$ is given by $\hat{\xi}_k=\xi(\beta_k)$ for a function $\xi\in\mathcal{H}_K$.  Then the minimum norm solution $\hat{g}^*$ that satisfies $W\hat{g}=\hat{\xi}$ is given by $\hat{g}^*=W' z$ where $WW' z= \hat{\xi}$ \citep*{luenberger69}, so that applying basic properties of the singular value decomposition yields
\begin{align} \label{appsol}
\hat{g}^*=\sum_{j=1}^{mq}\frac{\hat{\xi}' \bar{u}_j}{s_j}\bar{v}_j.
\end{align}
The components of the synthesized minimum norm control $\hat{u}^*=(\hat{u}_1^*,\ldots,\hat{u}_m^*)'$ are therefore given by
\begin{align} \label{syncont}
\hat{u}_k^*=\sum_{j=1}^{q}\frac{\hat{\xi}' \bar{u}_{k+m(j-1)}}{s_{k+m(j-1)}}\bar{v}_{k+m(j-1)}.
\end{align}
\begin{rem} The time and parameter discretizations $N$ and $P$ must be chosen such that $nP\leq mN$, so that the pair $(W,\hat{\xi}\,)$ represents an underdetermined system and therefore a minimum norm and not a least squares approximation problem.  The number $q$ of eigenfunctions used in the approximation is limited by $q\leq P$.
\end{rem}

\section{Examples and Simulations} \label{sec_example}
We present examples that illustrate the performance of our method.  We first consider the control of an ensemble of harmonic oscillators, which are frequently used as approximations for periodic phenomena in wide-ranging applications from engineering to physics, with parameter uncertainty as well as additive noise \citep*{bartlett02,mirrahimi04}.  Then, we demonstrate the control of an uncertain time-varying stochastic system, as well as a three-dimensional system of significance to quantum transport \citep*{stefanatos11}.

\begin{figure} 
\includegraphics[width=1.05\columnwidth]{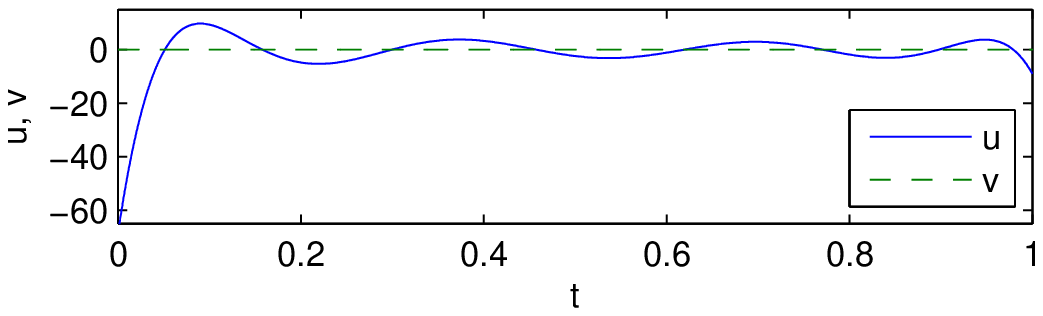} \put(-240,5){(a)}

\includegraphics[width=1.05\columnwidth]{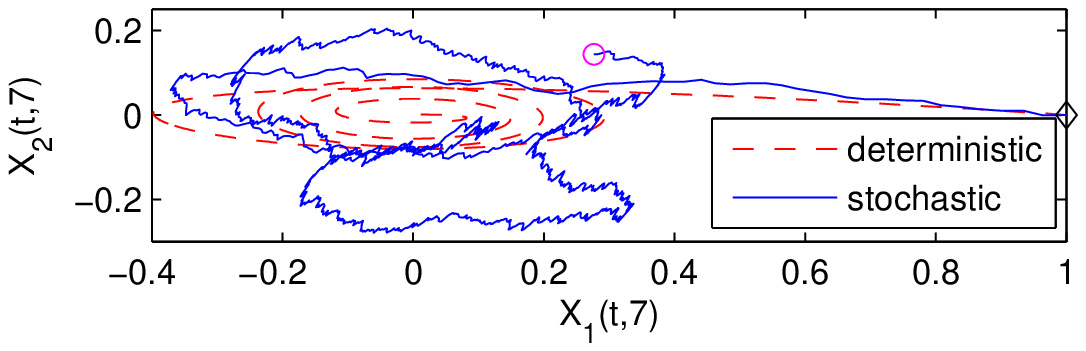} \put(-240,5){(b)}

\includegraphics[width=1.05\columnwidth]{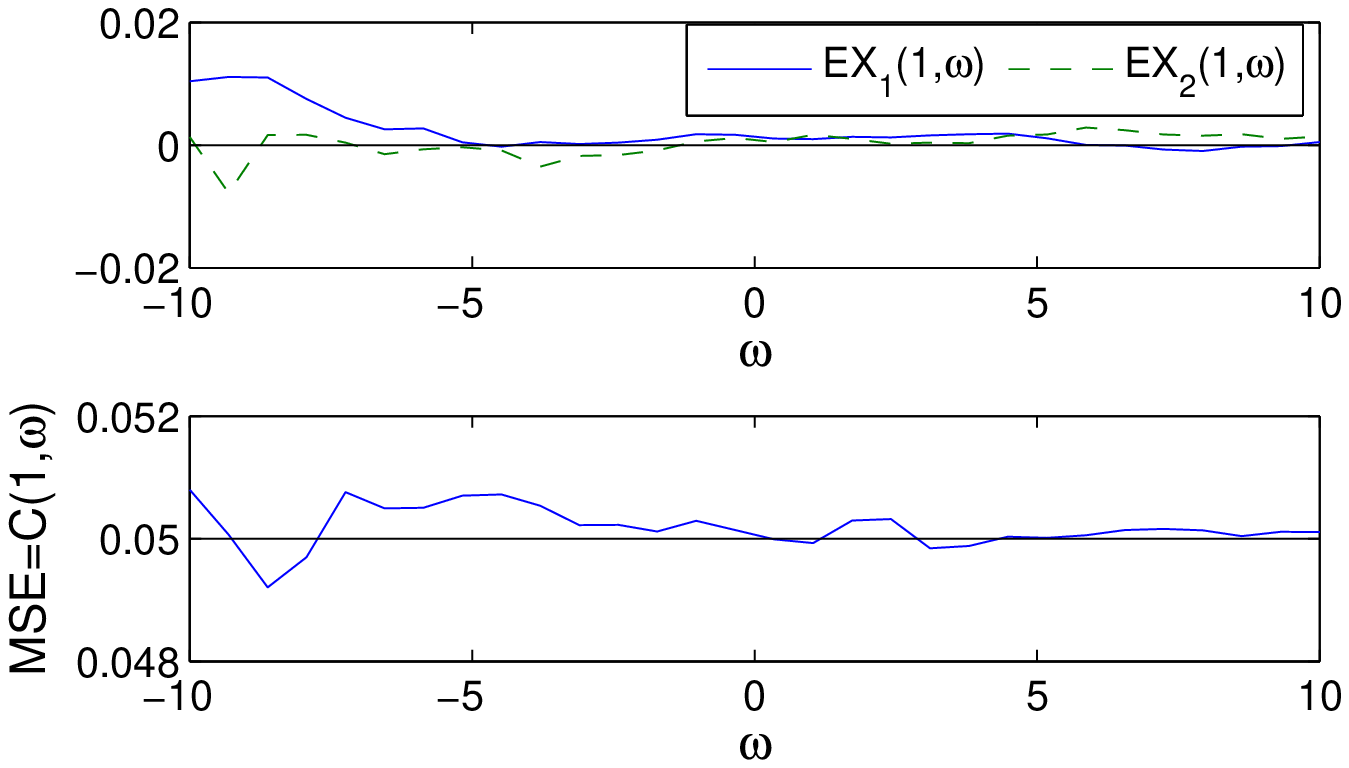} \put(-240,10){(d)} \put(-240,75){(c)}
\caption{Example \ref{example_ex1}, where $N=40001$ points on $[0,T]$ and $P=21$ samples of $\omega\in[-10,10]$ are used. (a) The optimal ensemble control. (b) A stochastic sample trajectory for $\omega=7$ using the Euler-Maruyama method with step $h=5\times 10^{-4}$, and the noise-free path. (c) The expectation and (d) MSE of the terminal state for 400 sample paths are close to the theoretical values of $(0,0)'$ and $0.05$, respectively.}
\label{fig_BM1}
\end{figure}

\begin{exmp} \label{example_ex1}
Consider a two-dimensional ensemble system with additive Gaussian process $\rd W$, given by
\begin{align}
\rd X(t,\omega)=A(t,\omega)X(t,\omega)\rd t+Bu(t)\rd t+G\rd W, \label{ex1}
\end{align}
$$
A(t,\omega)=\left[\begin{array}{cc} 0 & -\omega \\ \omega & 0 \\\end{array}\right], \quad B=\left[\begin{array}{cc} 1 & 0 \\ 0 & 1 \\\end{array}\right], \quad G=\left[\begin{array}{cc} 0.1 \\ 0.2 \\\end{array}\right],
$$
with frequency dispersion $\omega\in [-10,10]$.  We wish to steer this harmonic oscillator using a control $u(t)\in\mathbb{R}^2$ from initial state $X_0=(1,0)'$ to target state $X_F=(0,0)'$ at terminal time $T=1$.  The optimal ensemble control satisfies (\ref{critdet}), which is solved using the method in Section \ref{sec_operator}.  In order to avoid numerical conditioning errors, we choose $q=5$ in (\ref{appsol}) such that the largest and smallest singular values used satisfy $s_1/s_{mq}<10^4$ \citep*{zlotnik12acc}.  It is sufficient to  sample $\omega$ at 21 equidistant values on the interval $[-10,10]$ to obtain this number of eigenfunctions.  We use a discretization of 40001 points over the time interval $[0,1]$.  It is straightforward to compute the MSE at the terminal state using (\ref{indebm}) as $\mathrm{tr}\, C(T,\omega) = TG'G = 0.05$, which is invariant with respect to $\omega$.  The optimal control is shown in Figure \ref{fig_BM1}(a).  For the deterministic system, with $G=0$, it yields a terminal state error of less than $10^{-3}$ for all $\omega\in[-10,10]$.  The stochastic system is simulated using the Euler-Maruyama method \citep*{kloeden99} with sample trajectories shown in Figure \ref{fig_BM1}(b), and the statistics at the final state are shown to closely approximate the theoretical result in Figure \ref{fig_BM1}(c,d).
\end{exmp}

\begin{figure} 
\includegraphics[width=1.05\columnwidth]{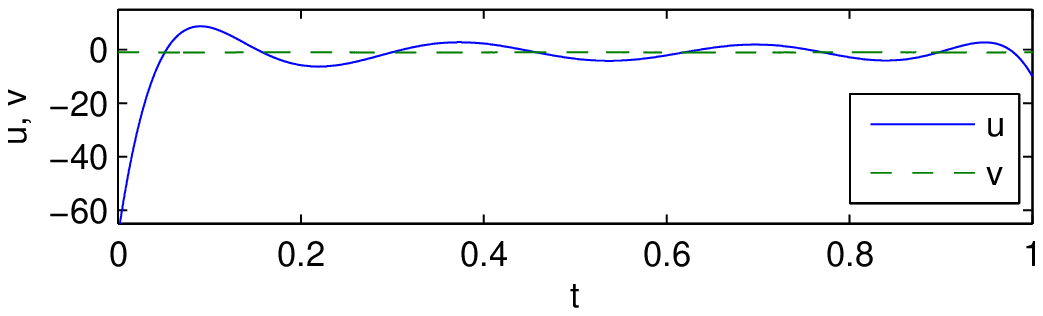} \put(-240,10){(a)}

\includegraphics[width=1.05\columnwidth]{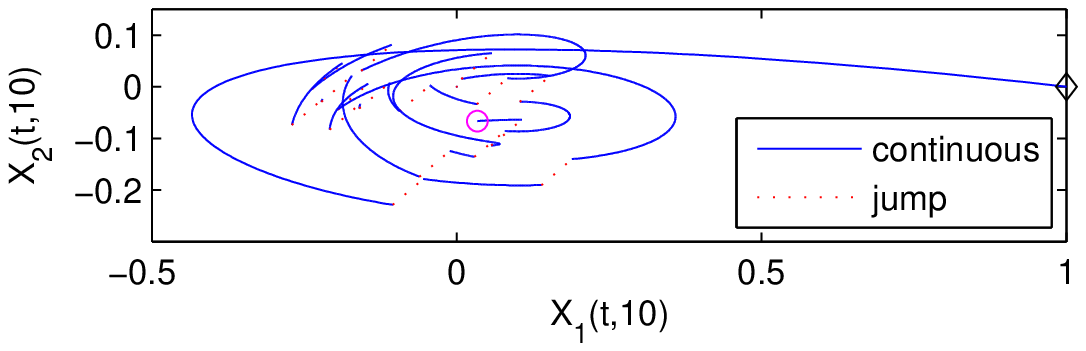} \put(-240,5){(b)}

\includegraphics[width=1.05\columnwidth]{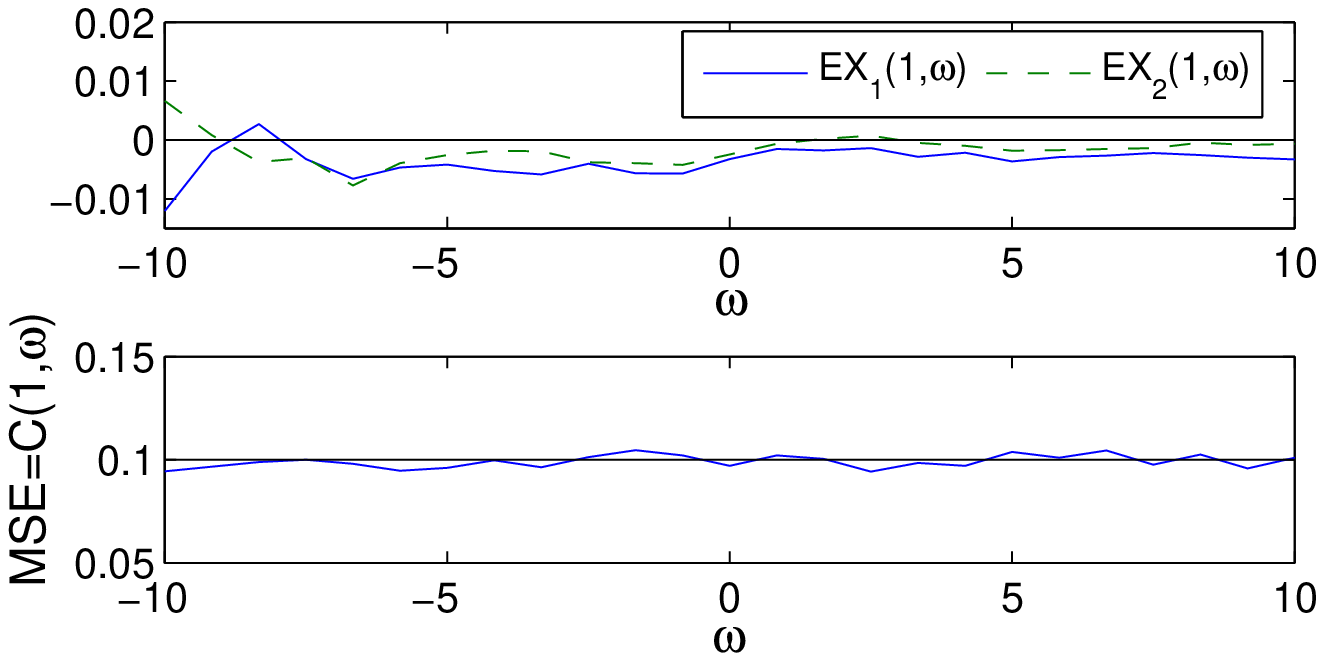} \put(-240,5){(d)} \put(-240,65){(c)}
\caption{Example \ref{example_ex2}, where $N=40001$ time nodes on $[0,T]$ and $P=21$ samples of $\omega$ on $[-10,10]$ are used. (a) The optimal ensemble control. (b) A stochastic sample trajectory for $\omega=-10$. (c) The expectation and (d) the MSE of the terminal state for 400 sample paths are close to the theoretical values of $(0,0)'$ and $0.1$, respectively.}
\label{fig_PP1}
\end{figure}

\begin{exmp} \label{example_ex2}
Consider the system
\begin{align}  \label{ex1pc}
\rd X(t,\omega)=A(t,\omega)X(t,\omega)\rd t+Bu(t)\rd t+G\rd N,
\end{align}
where $A$, $B$, $X_0$, and $X_F$ are as in Example \ref{example_ex1}, $\omega\in[-10,10]$, $G=(0.05,0.05)'$, $T=1$, and $\rd N$ is a scalar Poisson jump process with jump rate $\lambda=20$.  The optimal control is shown in Figure \ref{fig_PP1}(a). It is straightforward to compute the MSE at the terminal state using (\ref{indepp}) as $\mathrm{tr}\, C(T,\omega) = T\lambda G'G = 0.1$.  The Poisson jump system is simulated using a 4th order Runge-Kutta method and pseudo-randomly generated jump times, with sample trajectories shown in Figure \ref{fig_PP1}(d), and the statistics at the final state are shown to closely approximate the theoretical result in Figure \ref{fig_PP1}(b,c).
\end{exmp}

\begin{figure} 
\includegraphics[width=1.05\columnwidth]{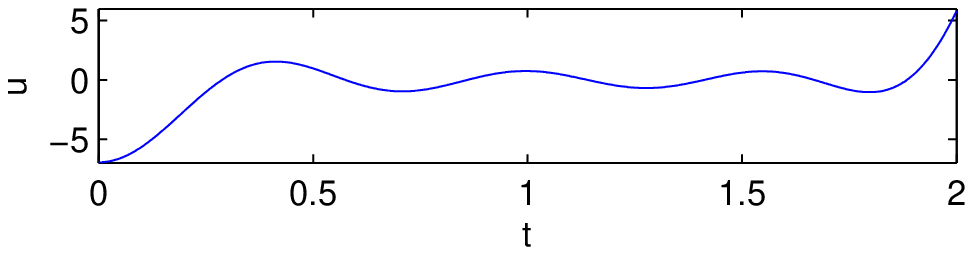} \put(-240,10){(a)}

\includegraphics[width=1.05\columnwidth]{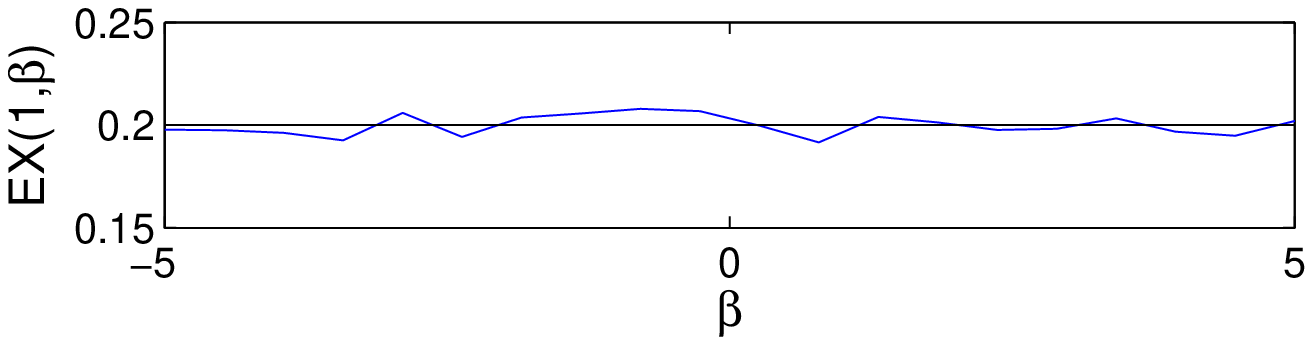} \put(-240,5){(b)}

\includegraphics[width=1.05\columnwidth]{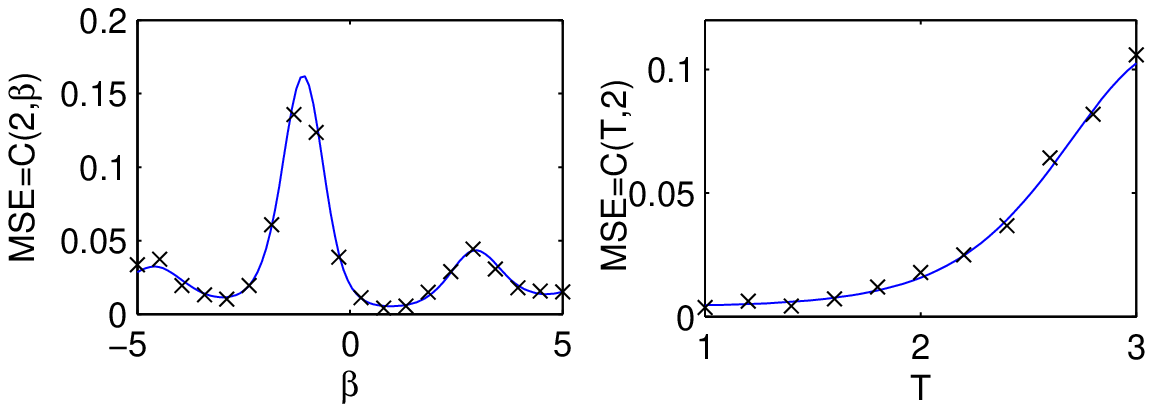} \put(-240,10){(c)}

\includegraphics[width=\columnwidth]{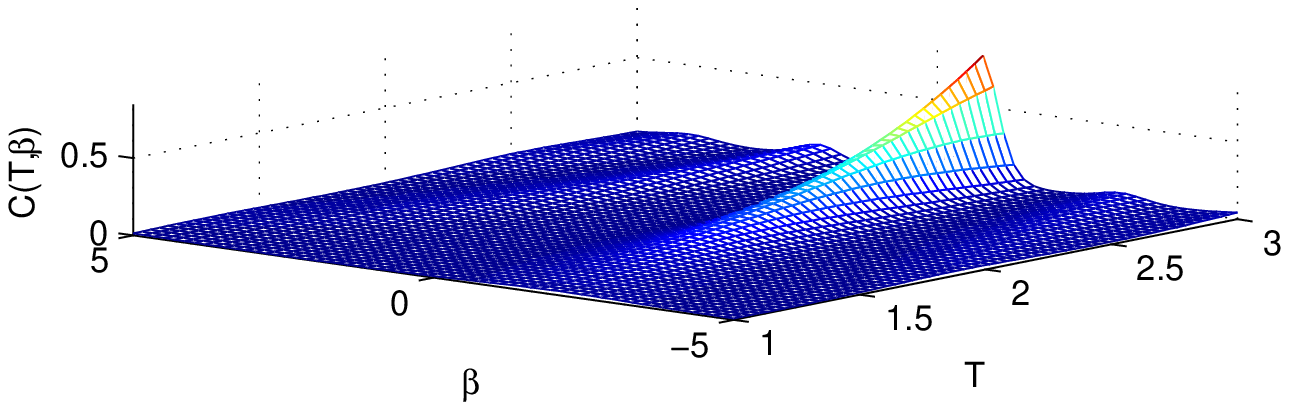} \put(-240,10){(d)}

\caption{Example \ref{example_ex3}, where $N=20001$ nodes on $[0,T]$ and $P=101$ samples on $\beta\in[-5,5]$ are used. (a) The optimal ensemble control, with $q=9$.  (b) The mean in the terminal state. (c) The MSE of the terminal state for $\beta\in[-5,5]$ when $T=2$, and for $T\in[1,3]$ when $\beta=2$.  The statistic computed using 400 sample paths with time step $h=10^{-3}$ ($\times$) is compared to the theory (line). (d) Surface plot of the theoretical MSE $C(T,\beta)$ for $(T,\beta)\in[1,3]\times[-5,5]$.}
\label{fig_TV1}
\end{figure}

\begin{exmp} \label{example_ex3} Consider a scalar time-varying system
\begin{align}
\rd x(t,\beta)=-\sin(\beta t)x(t,\beta)\rd t+u(t)\rd t+\rd W, \label{ex2a}
\end{align}
with $[\beta_1 ,\beta_2]=[-5,5]$, $T=1$, initial state $X_0=1$ and target state $X_F=0.2$.
The transition matrix for the system $\dot{x}(t,\beta)=-\sin(\beta t)x(t,\beta)$ is given by $\Phi(t,t_0,\beta)=\exp\int_{t_0}^t(-\sin(\beta\sigma)) \rd \sigma=\exp[\frac{1}{\beta}\cos(\beta t)-\frac{1}{\beta}\cos(\beta t_0)]$.
The optimal control satisfying (\ref{critdet}) is shown in Figure \ref{fig_TV1}(a). The expected MSE in the terminal state is obtained using (\ref{indebm}) as
$C(T,\beta)=\int_0^T\exp[\frac{2}{\beta}\cos(\beta t)-\frac{2}{\beta}\cos(\beta t_0)] \rd\tau$.  The expectation of the final state is close to the target as shown in Figure \ref{fig_TV1}(b).  For terminal time  $T=2$ the MSE for values of $\beta\in[-5,5]$ closely matches $C(2,\beta)$, and for $\beta=2$ the MSE for $T\in[1,3]$ closely match $C(T,2)$, as shown in Figure \ref{fig_TV1}(c,d).  Each statistic is estimated using 100 sample paths.  A surface plot of $C(T,\beta)$ on the region $(T,\beta)\in[1,3]\times[-5,5]$ is shown in Figure \ref{fig_TV1}(e).
\end{exmp}


\begin{figure} 
\includegraphics[width=1.05\columnwidth]{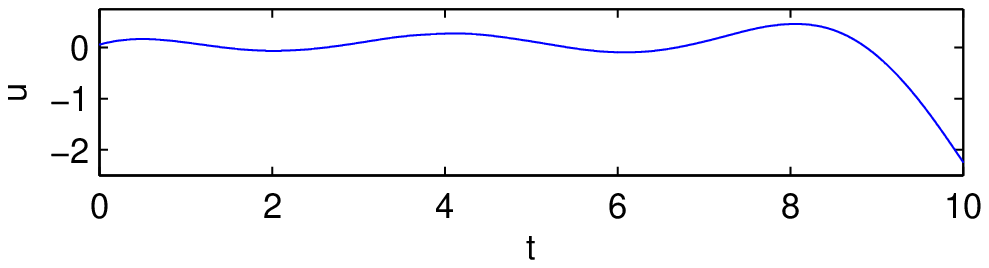} \put(-240,10){(a)}

\includegraphics[width=1.05\columnwidth]{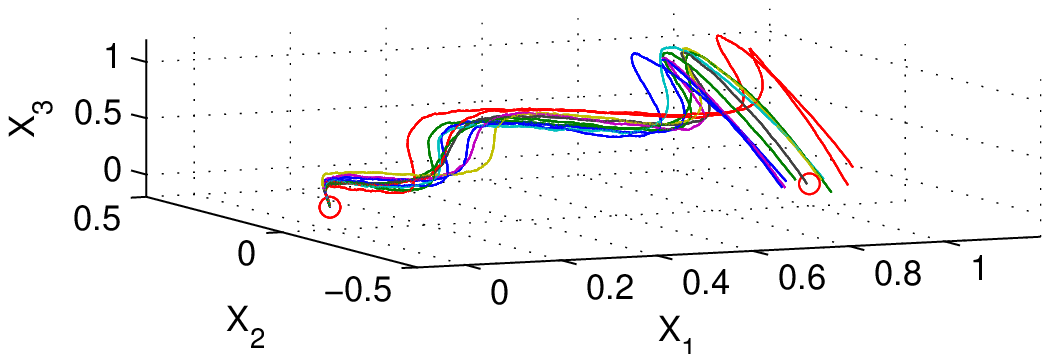} \put(-240,5){(b)}

\includegraphics[width=1.05\columnwidth]{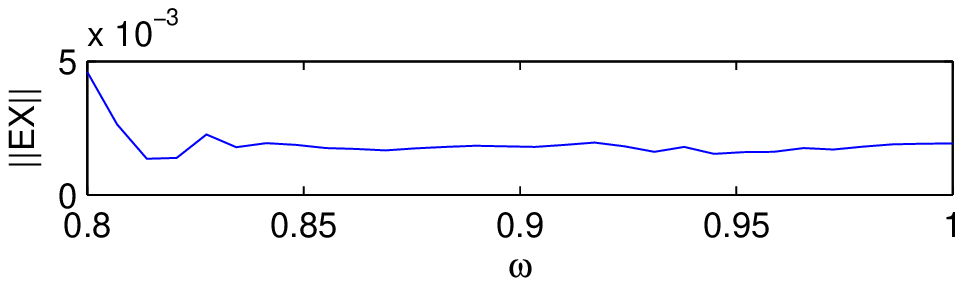} \put(-240,10){(c)}

\includegraphics[width=1.05\columnwidth]{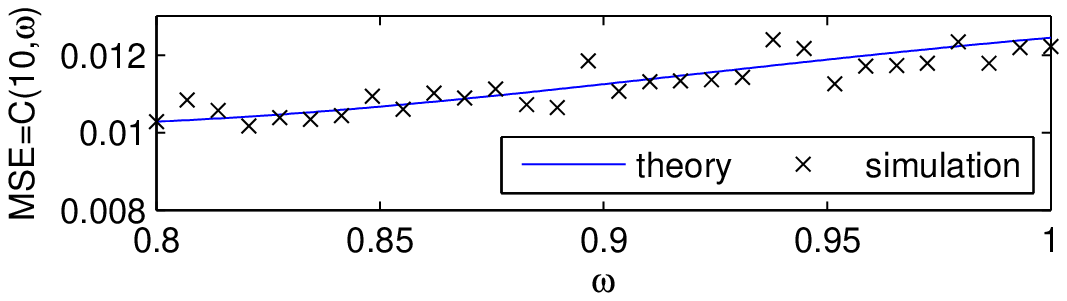} \put(-240,0){(d)}

\caption{Example \ref{example_ex4}, where $N=40001$ time nodes on $[0,20]$ and $P=101$ samples of $\omega$ on $[0.8,1]$ are used. (a) The optimal ensemble control.  (b) Sample trajectories for $\omega=0.9$. (c) The norm of the mean in the terminal state.   (d) The MSE of the terminal state. The statistic computed using 1000 sample paths with a time step $h=10^{-3}$.}
\label{fig_QT1}
\end{figure}

\begin{exmp} \label{example_ex4} A quantum ensemble transport system \citep*{stefanatos11} with noisy input is given by
\begin{align} \label{ex4}
\rd X=\left[A(t,\omega) X+ B u\right]\rd t + G \rd W,
\end{align}
$$
A(t,\omega)=\left[\begin{array}{ccc} 0 & 1 & 0 \\ -\omega^2 & 0 & \omega^2  \\ 0 & 0 & 0 \end{array}\right], \quad B=\left[\begin{array}{c} 0 \\ 0 \\ 1 \end{array}\right], \quad G=\left[\begin{array}{c} 0 \\ 0 \\ \sigma \end{array}\right],
$$
where $\rd W$ is a Gaussian process, $\omega\in[0.8,1]$, $\sigma=0.02$, $T=10$, and $X_0=(0,0,1)'$ and $X_F=(0,0,0)'$ are initial and target states, respectively. The optimal control satisfying (\ref{critdet}) is shown in Figure \ref{fig_QT1}(a), and sample trajectories are given in Figure \ref{fig_QT1}(b).  The statistics at the final state are computed from 1000 sample paths generated using a strong 1.5 order stochastic integration method \citep*{kloeden99}, and are shown in Figure \ref{fig_QT1}(c,d) to closely match the theory.  For terminal time  $T=10$, the MSE for all $\omega\in[0.8,1]$ closely matches $C(10,\omega)$, which is computed numerically by integrating (\ref{indebm}), and the mean terminal state is also near the target.
\end{exmp}

\section{Conclusion}  \label{sec_conc}
We have derived the optimal open-loop control that guides a family of independent, structurally identical, finite-dimensional stochastic linear systems with variation in system parameters, between initial and target states in function space on the parameter domain. Ensemble control has been extended to stochastic systems with additive Brownian noise and Poisson jump processes. It was shown that the same control minimizes both  the error in the mean and the mean square error in the terminal state with respect to the target in function space. The optimal control was obtained by solving the Fredholm integral equation associated with the system dynamics, and was approximated by using the singular value decomposition of a matrix that approximates the action of the Fredholm operator.  We used Monte Carlo stochastic integration to evaluate the statistical properties of the controlled example systems, and the results closely followed our theory with respect to performance objectives. Our work has immediate practical applications to the control of dynamical systems with additive noise and parameter dispersion, which are of interest in various areas such as NMR, MRI, and neuroscience, and also makes a fundamental contribution to stochastic control theory. The novel concepts we explored lead to a rich variety of new stochastic control problems, involving uncertainty in the system parameters and including optimization with respect to the state, control, and time horizon, for which standard methods, such as maximum principle and stochastic dynamic programming, may not be effective to obtain optimal controls.

\bibliographystyle{apa1}        

\footnotesize
\bibliography{stochastic_ref}

\end{document}